\input amstex
\documentstyle{amsppt}
\magnification=\magstep1
\tolerance=2000
\TagsOnRight
\def\gp#1{\langle #1 \rangle}
\def\m1{^{-1}}

\topmatter
\title
On the unitary subgroup of modular group algebras
\endtitle
\author
Victor Bovdi, Tibor  Rozgonyi
\footnote {Research supported by the Hungarian National Foundation
for Scientific Research No. 1903.}
\endauthor
\abstract
It this note we  investigate the structure of the group of $\sigma$-unitary
units in some noncommutative  modular group algebras $KG$, where $\sigma$ is
a non-classical ring involution of $KG$.
\endabstract
\endtopmatter

\document

Let $K$ be a field of characteristic $p$ and  $KG$   be  the
group algebra of the finite $p$-group $G$ over $K$. The  group
of units of $KG$ is denoted by $U(KG)$ and its subgroup
$$
V(KG) = \{\sum_{g\in G}\alpha_g g \in U(KG) \mid  \sum_{g\in G} \alpha_g = 1\}
$$
is called the group of normalized units.
Let $\sigma$ be an anti-automorphism of order two of the group $G$.
If $x =\sum_{g\in G}\alpha_gg$ is an element of the algebra $KG$,
then $u^\sigma$ denote the element $\sum_{g\in G}\alpha_gg^\sigma$.
The map $u\longmapsto u^\sigma$ is called $\sigma$-involution  and the
following relations hold:

\item{1.} $(u+v)^\sigma=u^\sigma+v^\sigma$;

\item{2.} $(uv)^\sigma=v^\sigma u^\sigma$;

\item{3.} $(u^\sigma)^\sigma=u$,
for all $u,v \in KG$.

We give example such involutions.

\item{1.} Let $\sigma (g)=g^{-1}$ for all $g\in G$. It is obvious that
$x\longmapsto x^\sigma$ is a classical involution of $KG$ and $x^\sigma$
is denoted by $x^*$.
\item{2.} Let $G$ be a finite 2-group and $C$ be the center of $G$ such that
$G/C\cong C_2\times C_2$ and the commutator
subgroup $G'=\gp{e\mid e^2=1}$ is of order two. Then the mapping
$\odot:G\longmapsto G$, defined by
$$
g^\odot=\cases $g,$  &\text{ if $g\in C$,}\\
                $ge,$ &\text{ if $g\notin C$}
         \endcases
$$
is an anti-automorphism of order two. If $x=\sum_{g\in G}\alpha_g g \in
KG$, then $u\longmapsto u^\odot=\sum_{g\in G}\alpha_g g^\odot$ is an
involution.

An element $u \in V(KG)$ is  called $\sigma$-unitary  if  $u^{-1}  =
u^\sigma$. The set of all unitary elements of the  group  $V(KG)$
is a subgroup. We shall denote this  subgroup  by  $V_\sigma (KG)$
and refer to $V_\sigma (KG)$ as the unitary subgroup  of  $V(KG)$.
It is clear that G is a subgroup of $V_\sigma (KG)$.

A.A.  Bovdi  and  A.A.  Szak\'acs  in \cite{1}  describe  the  group
$V_*(KG)$ for an arbitrary finite abelian $p$-group G and a finite field
K of characteristic $p$. A.A.  Szak\'acs  calculated the unitary subgroup
of commutative group algebras in \cite{4-6}. In case of a non-classical
involution $\sigma$ nothing is known of $V_\sigma (KG)$.

It is an open and important problem to investigate $V_\sigma (KG)$ in
the non-abelian case.   The  following  results  we obtain may be
considered the first ones in  this direction.
\proclaim
{Theorem 1} {\it Let $G$ be  a  finite  2-group  which  contains  an
abelian normal subgroup $A$ of index two. Suppose that there
exists an element $b \in G\setminus A$ of order $4$ such that
$b^{-1}ab = a^{-1}$ for all $a \in A$. Then  the  unitary
subgroup $V_*(\Bbb Z_2G)$ is the semidirect product of $G$ and  a
normal subgroup $H$.

The subgroup $H$ is  the  semidirect  product  of  the  normal
elementary abelian 2-group
$$
W = \{1 + ( 1 + b^2 )zb \mid z\in \Bbb Z_2A\}
$$
and the abelian subgroup $L$, where $V_*(\Bbb Z_2A) = A\times L$.

The  abelian  group   $W$   is   the   direct   product   of
$\frac12|A|$ copies of the additive group of field $\Bbb Z_2$.}
\endproclaim
{\bf Proof.\/} It is well known \cite{1} that
$$
V_*(\Bbb Z_2A) = A\times L.\tag1
$$
Suppose now that $G = \gp{A,b}$,  and
$b^{-1}ab = a^{-1}$ for  all  $a\in A$.  Then  every  element  of
$V_*(ZG)$  has  a unique representation of form
$$
                 x = x_1 + x_2b \hskip10pt       (x_i \in \Bbb Z_2A).
$$
Let $\chi(x_i)$ be the sum of the  coefficients  of  the  element
$x_i$.

Let $J(A)$ denote the ideal  of  the  ring  $\Bbb Z_2G$  generated  by the
elements   of   form   $h-1$   with    $h\in A$.  Since
$\Bbb Z_2G/J(A)\cong \Bbb Z_2(G/A)$ and quotient group $G/A$  is  of   order
$2$, it   follows   that   element $$
\overline{x}=\chi(x_1)+\chi(x_2)\overline{b} \hskip10pt (\overline{b}=bA)
$$
is trivial. Hence, one of the elements $\chi(x_1)$ or
$\chi(x_2)$  equals  $1$  and  the  other  is  zero. Since
$G\subseteq V_*(\Bbb Z_2G)$, element  $xb  \in V_*(\Bbb Z_2G)$ and $x$ or $xb$
have the form
$$
x_i(1  +  x_i^{-1}x_jb),\hskip10pt  (i,j  \in \{1,2\}, i\not=j).
$$
Indeed, if $\chi(x_i) = 1$, then $x_i$ is a unit.

Suppose that $x = x_1(1 + x_2b) \in  V_*(\Bbb Z_2G)$  and  $x_1  \in
V(\Bbb Z_2A)$  and $\chi(x_2) = 0$. Then
$$
x^* = (1 + b^{-1}x_2^*)x_1^* \hskip8pt and \hskip8pt xx^* = 1.
$$
It follows that
$$
\cases
\text{$x_1x_1^*(1 + x_2x_2^*) = 1$,}\\
\text{$x_2(1+b^2)=0$.} \\
\endcases     \tag2
$$
Then by Proposition of \cite{2}  we have
$$
x_2 = (1 + b^2)z \hskip10pt (z \in  \Bbb Z_2A).
$$
It is obvious that
$$
x_2x_2^* = 2(1 + b^2)zz^* = 0.
$$
From (2) we have  $x_1x_1^*= 1$ and element $x_1$ is unitary. Therefore
$x_1\in  V_*(\Bbb Z_2A)$.

Let $g_1,\ldots ,g_t$ be the representatives  of  the  distinct cosets of
$A$ modulo $\gp{b^2}$ and
$$
W_{g_i}= \{1 + (1+b^2)\alpha g_ib
\mid\alpha  \in  \Bbb Z_2\}.
$$
It is clear that $W_{g_i}$ is a group of order $2$ and $W =\prod_{i=1}^z W_{g_i}$
is a direct product. If $w_i = 1 + (1 + b^2)g_ib$ and
$bg_ib^{-1}=g_i^{-1}=g_jb^k$ \hskip8pt ($k\geq 2$) then
$$
bw_ib^{-1} = 1 +(1 + b^2)g_jb \in W \tag3
$$
and for every $x_1 \in V_*(\Bbb Z_2A)$
$$
x_1w_ix_1^{-1} =1+(1+b^2)x_1^2g_ib \in  W, \tag4
$$
because
$$
bx_1^{-1}= b{x_1}^* = x_1b    \tag5
$$
Thus $W$ is a normal subgroup of $V_*(\Bbb Z_2G)$ and $W \cap L=1$, where
$V_*(\Bbb Z_2A) = A\times L$.

Let $H$ be the subgroup generated by $W$ and $L$. By (4) it  follows
that $H$ is a semidirect product of normal subgroup $W$ and the
subgroup $L$.

Therefore we proved that $V_*(\Bbb Z_2G)$  is  generated  by  $G$
and  $H$. From (3), (4) and (5) it follows that H is a normal  subgroup  of
$V_*(\Bbb Z_2G)$.

\proclaim
{Theorem 2} {\it Let $C$ be the center of a finite 2-group $G$ and

\item{1.} $G/C$ is a direct product of two groups $\gp{aC}$ and $\gp{bC}$ of order two;
\item{2.} the commutator subgroup $G'=\gp{e}$ of group $G$ is order $2$;

Then the unitary subgroup
$$
V_\odot (\Bbb Z_2G)=G\times T\times W,
$$
where $W=\{1+x_1a+x_2b+x_3ab \mid x_i\in \Bbb Z_2C(1+e)\}$ is a central
elementary abelian 2-group, the subgroup $V(\Bbb Z_2C)[2]$ of
all elements of order $2$ in $V(\Bbb Z_2C)$ is a direct product of group $T$
and the subgroup $C[2]$ of all elements of order $2$ in $C$.}
\endproclaim
{\bf Proof.\/} Let $a,b\in G$ and $[a,b]=e$. Then $C=C_G(a,b)$ is the
center of $G$ and
$$
G=C\cup Ca \cup Cb \cup Cab.
$$
Every element of $V_\odot (\Bbb Z_2G)$ can be written in the form
$$
x=x_0+x_1a+x_2b+x_3ab \hskip10pt (x_i\in \Bbb Z_2C)
$$
and elements $x_i$ are central in $KG$. Then
$$
x^\odot=x_0+x_1ae+x_2be+x_3abe \hskip10pt (x_i\in \Bbb Z_2C)
$$
and a simple calculation immediately shows that $x\in V_\odot(\Bbb Z_2G)$
if and only if
$$
\cases
\text{$x_0^2+x_1^2a^2e+x_2^2b^2e+x_3^2a^2b^2=1$,}\\
\text{$(x_0x_1+x_2x_3b^2)(1+e)=0$,}\\
\text{$(x_0x_2+x_1x_3a^2)(1+e)=0$,}\\
\text{$(x_0x_3+x_1x_2)(1+e)=0$.}\\
\endcases \tag6
$$
Since $\chi(x_i)$ equals 1 or 0 and $G\subseteq V_\odot (\Bbb Z_2G)$ we
may assume that $\chi(x_0)=1$. By Propositions 2.7 \cite{2} from (6),
we obtain
$$
x_0x_1+x_2x_3b^2=(1+e)r_1,
$$
$$
x_0x_2+x_1x_3a^2=(1+e)r_2,
$$
$$
x_0x_3+x_1x_2=(1+e)r_3,
$$
for some $r_i\in \Bbb Z_2C$.

Multiplying the first equality  by $x_1a^2$ and the second equality
by $x_2b^2$ we have
$$
x_0(x_1^2a^2+x_2^2b^2)=(1+e)r_4,
$$
$$
x_0(x_1^2+x_3^2b^2)=(1+e)r_5,
$$
$$
x_0(x_2^2+x_3^2a^2)=(1+e)r_6,
$$
for some $r_i \in \Bbb Z_2C$.
Furthermore,  Proposition 2.6 \cite{2} assures that $\chi$ is a ring
homomorphism. Because $\chi (x_0)=1$, then we have
$$
\cases
\text{$\chi (x_1)+\chi (x_2)=0,$}\\
\text{$\chi (x_1)+\chi (x_3)=0,$}\\
\text{$\chi (x_2)+\chi (x_3)=0.$}\\
\endcases
$$
Clearly, if $x\in V(\Bbb Z_2G)$ then $\chi (x)=\chi (x_0)+\chi (x_1)+\chi
(x_2)+ \chi (x_3)=1$ and therefore it follows that $\chi
(x_1)=\chi (x_2)= \chi (x_3)=0$.

Since $x_0$ is a unit, $x$ can be written as
$$
x=x_0(1+y_1a+y_2b+y_3ab),
$$
where $\chi (y_i)=0$ and $y_i\in \Bbb Z_2C$. Then
$$
xx^\odot=x_0^2(1+y_1a+y_2b+y_3ab)(1+y_1ae+y_2be+y_3abe)=1
$$
if and only if
$$
\cases
\text{$x_0^2(1+y_1^2a^2e+y_2^2b^2+y_3^2a^2b^2)=1$,}\\
\text{$y_1(1+e)+y_2y_3b^2(1+e)=0,$}\\
\text{$y_2(1+e)+y_1y_3a^2(1+e)=0,$}\\
\text{$y_3(1+e)+y_1y_2(1+e)=0.   $}\\
\endcases \tag7
$$
Then from second and fourth equalitis we obtain
$$
y_1(1+e)=y_1y_2^2b^2(1+e)
$$
and
$$
y_1(1+e)(1+y_2^2b^2)=0.
$$
Since $\chi (y_2)=0$, element $1+y_2^2b^2$ is a unit and from the
equality above it follows that $y_1(1+e)=0$. Clearly, (7) implies
$y_2(1+e)=y_3(1+e)=0$. By virtue of Proposition 2.7 \cite{2}
$y_i=(1+e)u_i$ ($i=1,2,3$) and $u_i\in \Bbb Z_2C$.

Clearly, $y_i^2=0$ ($i=1,2,3$) and this together with (7) gives
$x_0^2=1$. It is obvious that
$$
y=1+y_1a+y_2b+y_3ab \hskip10pt (y_i\in (1+e)\Bbb Z_2C)
$$
is a central unit of order $2$.

It is well known that
$$
V(\Bbb Z_2C)[2]=C[2]\times T
$$
and $T\cap G=1$.

Clearly, we proved that $V_\odot (\Bbb Z_2G)$ is generated by $G$ and two
central subgroups:
$$
V(\Bbb Z_2C)[2]=\{x_0\in V(\Bbb Z_2G) \mid x_0^2=1\},
$$
$$
W=\{1+y_1a+y_2b+y_3ab \mid y_i\in (1+e)\Bbb Z_2C \}.
$$
Therefore
$$
V_\odot (\Bbb Z_2G)=G\times T\times W.
$$

\address
Department of Mathematics, Bessenyei Gy\"orgy Teachers' Training College,
Ny\'\i regyh\'aza, Hungary
\endaddress
\enddocument